\documentclass[12pt,a4paper,twoside]{article}
\usepackage[UKenglish]{babel}
\usepackage[T1]{fontenc}
\usepackage[utf8x]{inputenc}
\usepackage{latexsym,amsfonts,amsmath,amsthm,amssymb,mathrsfs}
\usepackage{fullpage}
\usepackage[table]{xcolor}
\usepackage{bbm}
\usepackage{paralist}
\usepackage{todonotes}
\usepackage{dsfont}
\usepackage{float}
\usepackage{caption}
\usepackage{xfrac}
\usepackage[normalem]{ulem}
\newcommand{\supp}{\operatorname{supp}}
\usepackage{relsize}
\usepackage{mathtools}
\usepackage{array}
\usepackage{makecell}
\usepackage{multirow}

\usepackage{cancel}
\usepackage{soul}
\usepackage{dsfont}
\usepackage[colorlinks=true, linkcolor=teal, citecolor=-teal,
pagebackref=true]{hyperref}

\usepackage{graphicx}

\usepackage{fourier}

\newcommand{\yq}{(\mathbf{y}_q)_{q \in \mathbb{N}} }

\newcommand{\by}{\mathbf{y}}
\newcommand{\byq}{\mathbf{y}_q}
\newcommand{\Aq}{A_q^{y}}
\newcommand{\Aqm}{\mathbf{A}_q^{\by}}
\newcommand{\Arm}{\mathbf{A}_r^{\by}}
\newcommand{\Ar}{A_r^{z}}
\newcommand{\overlap}{\lambda_1(\Aq \cap \Ar)}

\newcommand{\productm}{\left(\frac{\psi(q)\varphi(q)}{q}\frac{\psi(r)\varphi(r)}{r}\right)^m}

\DeclareMathOperator{\R}{\mathbb{R}}
\DeclareMathOperator{\N}{\mathbb{N}}

\DeclareMathOperator{\lcm}{lcm}
\DeclareMathOperator{\Log}{Log}

\def\eps{{\varepsilon}}

\newcommand{\calH}{\mathcal H}
\newcommand{\ba}{\mathbf a}
\newcommand{\NN}{\mathbb N}
\newcommand{\q}{\mathbf q}
\newcommand{\RR}{\mathbb R}

\newcommand{\x}{\mathbf x}
\newcommand{\X}{\mathbf X}
\newcommand{\y}{\mathbf y}
\newcommand{\ZZ}{\mathbb Z}

\newcommand{\measqm}{\left(\frac{\psi(q)\varphi(q)}{q}\right)^m}
\newcommand{\measr}{\frac{\psi(r)\varphi(r)}{r}}
\newcommand{\measrm}{\left(\frac{\psi(r)\varphi(r)}{r}\right)^m}

\DeclarePairedDelimiter{\abs}{\lvert}{\rvert}
\DeclarePairedDelimiter{\parens}{\lparen}{\rparen}
\DeclarePairedDelimiter{\set}{\lbrace}{\rbrace}

\newtheorem{thm}{Theorem}

\newtheorem{lem}{Lemma}[section]
\newtheorem{cor}[lem]{Corollary}

\newtheorem{prop}[lem]{Proposition}

\newtheorem{conj}{Conjecture}
\newtheorem{question}{Question}

\theoremstyle{remark}
\newtheorem*{rem}{Remark}

\allowdisplaybreaks

\title{The Duffin--Schaeffer conjecture with a moving target}
\author{Manuel Hauke \\ University of York \and Felipe A. Ram\'irez \\ Wesleyan University}

\date{}

\begin{document}

\maketitle

\begin{abstract}
We prove the inhomogeneous generalization of the Duffin--Schaeffer conjecture in dimension $m \geq 3$. That is, given $\y\in\RR^m$ and $\psi:\NN\to\RR_{\geq 0}$ such that $\sum (\varphi(q)\psi(q)/q)^m = \infty$, we show that for almost every $\x\in\RR^m$ there are infinitely many rational vectors $\ba/q$ such that $\abs{q\x - \ba - \y}<\psi(q)$ and such that each component of $\ba$ is coprime to $q$. This is an inhomogeneous extension of a homogeneous conjecture of Sprind\v{z}uk which was itself proved in 1990 by Pollington and Vaughan. In fact, our main result generalizes Pollington--Vaughan not only to the inhomogeneous case, but also to the setting of moving targets, where the inhomogeneous parameter $\y$ is free to vary with $q$. In contrast, we show by an explicit construction that the ($1$-dimensional) inhomogeneous Duffin--Schaeffer conjecture fails to hold with a moving target, implying that any successful attack on the one-dimensional problem must use the fact that the inhomogeneous parameter is constant. We also introduce new questions regarding moving targets.
\end{abstract}

\setcounter{tocdepth}{1}
\tableofcontents

\section{Introduction}

Given a function $\psi:\NN\to[0,\infty)$ and a real
number $y$, let $W'(y,\psi)$ be the set of real
numbers $x\in[0,1]$ for which the inequality
\begin{equation}\label{eq:1}
\abs*{qx - a - y} < \psi(q)
\end{equation}
has infinitely many solutions $(a,q)\in\ZZ\times\NN$ with
$\gcd(a,q)=1$. When there is no coprimality requirement, denote the set by $W(y,\psi)$.

Khintchine's Theorem~\cite[1924]{K24} says that
if $\psi$ is monotonic\footnote{The original statement of Khintchine needed the stronger assumption of $q\psi(q)$ being monotonically decreasing, but it is well-known that $\psi$ being non-increasing is sufficient.} and $\sum \psi(q)$ diverges, then $W(0,\psi)$
has full Lebesgue measure. This was generalized to the inhomogeneous
case, where $y\in\RR$ is arbitrary, by Szüsz~\cite[1958]{SzuszinhomKT}.

These results are also known to hold in higher dimensions.
Given $m\in\NN$ and a real vector $\y \in \RR^m$, we
define $W_m(\y,\psi)$ to be the set of 
$\x\in[0,1]^m$ for which the inequality
\begin{equation*}
\abs*{q\x - \ba - \y} < \psi(q)
\end{equation*}
has infinitely many solutions $(\ba,q)\in\ZZ^m\times\NN$, where
$\abs{\cdot}$ denotes maximum norm. With the additional coprimality condition
\begin{equation}\label{eq:coprimality}
    \gcd(a_i,q) = 1\qquad (i=1, \dots, m),
\end{equation}
the set is denoted $W_m'(\y,\psi)$.

In 1926 Khintchine~\cite{K26} proved the $m$-dimensional version of his theorem for $W_m(\mathbf{0},\psi)$, and in 1964~\cite{S64} Schmidt proved the $m$-dimensional inhomogeneous version for $W_m(\y,\psi)$. In these results, the series of interest is $\sum \psi(q)^m$.\\

The question of removing the monotonicity assumption in Khintchine's Theorem and its higher-dimensional and inhomogeneous companions is both classical and of recent interest. In their seminal work, Duffin and Schaeffer~\cite{duffinschaeffer} showed by counterexample that monotonicity cannot be completely removed from Khintchine's Theorem. In place of this, they conjectured a Khintchine-type theorem to hold when only considering reduced fractions. To be more precise, it was conjectured that for an arbitrary function $\psi:\NN\to\RR_{\geq 0}$, the set $W'(0,\psi)$ has full measure if and only if $\sum\varphi(q)\psi(q)/q$ diverges. The Duffin--Schaeffer
conjecture was proved in a breakthrough result by Koukoulopoulos and
Maynard~\cite{KMDS} in 2020. In a preceding work, Pollington and Vaughan~\cite{PV90} proved in 1990 an $m$-dimensional version for $m\geq 2$, confirming a conjecture of Sprind\v{z}uk~\cite{Sprindzuk}.\\ 

In this article, we consider inhomogeneous variants of the Duffin--Schaeffer conjecture, a topic that has gained quite a lot of interest in recent years. 
The inhomogeneous Duffin--Schaeffer conjecture postulates that for any real number $y$, the set $W'(y,\psi)$ has full measure if and only if $\sum\varphi(q)\psi(q)/q$ diverges. 
This was first explored in detail in~\cite{ds_counterex} where inhomogeneous versions of the Duffin--Schaeffer counterexample were provided. 
Although there has been partial progress when considering special approximation functions $\psi$~\cite{CT24} or certain inhomogeneous parameters~\cite{BHV24}, the conjecture itself remains wide open.\\

Our main result establishes the $m$-dimensional version of this conjecture for $m\geq3$, confirming the inhomogeneous version of Sprind\v{z}uk's conjecture in those dimensions. It can be seen as an inhomogeneous analogue to the work of Pollington and Vaughan~\cite{PV90}. In fact, we prove the statement more generally, with a moving target $\by = (\by_q)_{q\in\NN}$. Now $W_m(\y,\psi)$ and $W_m'(\y,\psi)$ are defined by
\begin{equation*}
\abs*{q\x - \ba - \y_q} < \psi(q)
\end{equation*}
being satisfied with infinitely many $q\in\NN$. A fixed inhomogeneous parameter $\y\in\RR^m$ can be interpreted as a constant sequence.

\begin{thm}\label{thm:movingdsc}
    Let $m\geq 3$. Let $\psi: \N \to [0,\infty)$ be an arbitrary function and $\mathbf{y} = \yq \subset \mathbb{R}^m$ be an arbitrary sequence.
    Denoting the $m$-dimensional Lebesgue measure by $\lambda_m$, we have
    \[\lambda_m\left(W_m
    '(\by,\psi)\right) =
\begin{cases}
  0 &\text{if } \sum\limits_{q \in \N} \left(\frac{\varphi(q)\psi(q)}{q}\right)^m < \infty,\\
    1 &\text{if } \sum\limits_{q \in \N} \left(\frac{\varphi(q)\psi(q)}{q}\right)^m = \infty.
\end{cases}
    \]
\end{thm}

\begin{rem}
This theorem contains the inhomogeneous problem as a special case, by setting $\y$ to be a constant sequence. 
\end{rem}

As a consequence of Theorem \ref{thm:movingdsc}, we obtain its analogue for $n$-by-$m$ systems of linear forms (see Theorem \ref{thm:bootstrapped}), solving the $m\geq 3$ cases of a conjecture appearing in~\cite[Conjecture~1]{AR_2023}, and leaving open the cases where $m =1$ and $m =2$. See Section~\ref{sec:remarks} for details on this conjecture and further results in the theory of systems of linear forms.

Since the statement of Theorem \ref{thm:movingdsc} allows for a moving target in the higher-dimensional inhomogeneous Duffin--Schaeffer conjecture, it becomes natural to ask whether the same is true in dimensions $m=1$ and $m=2$. Does the original one-dimensional inhomogeneous Duffin--Schaeffer conjecture allow for a moving target? Surprisingly, it turns out that for $m=1$, this general statement is false!

\begin{thm}\label{thm:weakcounterex}
  There exists a sequence $y = (y_q)_{q \in \N}$ of reals and corresponding
  $\psi:\NN\to[0,\infty)$ such that
  $\sum_{q \in \mathbb{N}}\frac{\varphi(q)\psi(q)}{q} = \infty$, yet $\lambda(W'(y, \psi))=0$.
\end{thm}

This leaves open only the $m=2$ case, which is discussed as a conjecture in Section~\ref{sec:remarks}, together with some other $2$-dimensional problems that have arisen in recent literature.

\subsection{Consequences, remarks and open questions}\label{sec:remarks}

\subsubsection*{The inhomogeneous Duffin--Schaeffer conjecture}

Although we do not see a way to adapt Theorem \ref{thm:weakcounterex} to the case where $\by$ is constant, the constructed counterexample can be seen as an argument for doubters of the inhomogeneous Duffin--Schaeffer conjecture. The statement shows that any possible proof of the inhomogeneous Duffin--Schaeffer conjecture has to take into account the fact that $\by$ is constant.

We remark that Theorem \ref{thm:weakcounterex} does not provide a counterexample to the so-called weak inhomogeneous Duffin--Schaeffer conjecture~\cite[Conjecture~1.22]{CT24}, not even in the case of a moving target. Therefore, it is still possible that the following statement holds true.

\begin{question}[Weak Duffin--Schaeffer conjecture with a moving target]\label{quest}
    Let $\psi: \N\to [0,\infty)$ be an arbitrary function with the property $\sum_{q \in \N} \frac{\psi(q)\varphi(q)}{q} = \infty$. Does one have for every sequence $y = (y_q)_{q \in \N}$ that $\lambda(W(y,\psi)) = 1$?
\end{question}

 Note that since $\liminf_{q \to \infty} \frac{\varphi(q)}{q} = 0$, this can be seen as a statement with a form of extra-divergence. Possibly, the use of the $\varphi$-factors can be justified by considering shifted forms of coprimality (that might depend on $y$), as it was done in~\cite{BHV24,CT24}.

\subsubsection*{On statements in dimension 2 and generalizations to linear forms}

Theorems \ref{thm:movingdsc} and \ref{thm:weakcounterex} leave open the question about the statement being true in dimension $m = 2$. We expect Theorem \ref{thm:movingdsc} also to hold in dimension $2$, also allowing a moving target.

\begin{conj}\label{conj:weakdscmoving}
     Let $\psi: \N \to [0,\infty)$ be an arbitrary function and $\mathbf{y} = \yq \subset \mathbb{R}^2$ be an arbitrary sequence.
    Then 
    \[\lambda_2\left(W_2'(\by,\psi)\right) =
\begin{cases}
  0 &\text{if } \sum\limits_{q \in \N} \left(\frac{\varphi(q)\psi(q)}{q}\right)^2< \infty,\\
    1 &\text{if } \sum\limits_{q \in \N} \left(\frac{\varphi(q)\psi(q)}{q}\right)^2 = \infty.
\end{cases}
    \]
\end{conj}

Possibly, the current methods allow one to show $\lambda_2(W_2(\by,\psi)) = 1$ when assuming some form of extra-divergence, e.g. assuming $\sum\limits_{q \in \N} \frac{\varphi(q)^{2+\eta}\psi(q)^2}{q^{2+\eta}} = \infty$ for some $\eta > 0$.
This would accord with comparable known results for the inhomogeneous Khintchine Theorem without monotonicity~\cite{allen2021independence,Yu}, where in dimension $m=2$, a similar form of extra-divergence is assumed. In~\cite{allen2021independence}, it is conjectured that the inhomogeneous Khintchine Theorem holds also without monotonicity in dimension $m=2$,  as it does for the homogeneous version---a classical result by Gallagher~\cite{G65}.\\

All the statements discussed so far have their generalizations to the setting of systems of $m$ linear forms in $n$ variables. Given $m,n\in\NN$, a sequence $\y = (\y_q)_{q \in \N}\subset \RR^m$ and a function $\psi:\ZZ^n \to [0,\infty)$, let $W_{n, m}(\y,\psi)$ be the set of 
$\X\in[0,1]^{nm}\subset \operatorname{Mat}_{n\times m}(\RR)$ for which the inequality
\begin{equation*}
\abs*{\q\X - \ba - \y_\q} < \psi(\q)
\end{equation*}
has infinitely many solutions $(\ba,\q)\in\ZZ^m\times\ZZ^n$. Let $W_{n,m}'(\y,\psi)$ denote the same set with the added condition that $\gcd(a_i,\q)=1$ for all $i=1, \dots, m$.

The analogue of Khintchine's Theorem in this context is known as the Khintchine--Groshev Theorem~\cite[1938]{Groshev}. The theory for systems of linear forms has developed in parallel with the classical ``simultaneous'' theory, but with a lag of a few years. Khintchine was followed by Khintchine--Groshev; a version of Khintchine's Theorem without monotonicity by Gallagher~\cite[1965]{G65} was followed by a linear forms version by Beresnevich--Velani~\cite[2010]{BV10}; Yu's $m\geq 3$ cases of Theorem \ref{thm:movingtarget}~\cite[2021]{Yu} were followed by the work of Allen and the second author~\cite[2023]{allen2021independence}; and the proofs of Pollington--Vaughan~\cite[1990]{PV90} and Koukoulopoulos--Maynard~\cite[2020]{KMDS} were recently followed with a proof of the linear forms version of the Duffin--Schaeffer conjecture~\cite{ramirez2023duffinschaeffer}. These results and others are summarized in Tables~\ref{KG} and~\ref{tableDS}. 

Recently, it has been shown that a large class of positive measure statements for limsup sets in $\RR^m$ can be bootstrapped to corresponding positive measure statements for limsup sets of $n$-by-$m$ systems of linear forms~\cite{ramirez2024metric}. So, for example, the linear forms results mentioned in the previous paragraph all follow from their earlier ``simultaneous'' versions. The bootstrapping results from~\cite{ramirez2024metric} can be used in conjunction with Theorem \ref{thm:movingdsc} to prove the following. 

\begin{thm}\label{thm:bootstrapped}
    Let $m,n\in\NN$ with $m\geq3$. Let $\psi: \ZZ^n \to [0,\infty)$ be an arbitrary function and $(\y_d)_{d=1}^\infty \subset \mathbb{R}^m$ be an arbitrary sequence and denote by $\y = (\y_\q)_{\q\in\ZZ^n}$ the targets $\y_\q = \y_{\gcd(\q)}$.
    Then 
    \[\lambda_{nm}\left(W_{n,m}'(\y,\psi)\right) =
\begin{cases}
  0 &\text{if } \sum\limits_{\q\in\ZZ^n\setminus \set{0}} \parens*{\frac{\varphi(\gcd(\q))\psi(\q)}{\gcd(\q)}}^m < \infty,\\
    1 &\text{if } \sum\limits_{\q\in\ZZ^n\setminus \set{0}} \parens*{\frac{\varphi(\gcd(\q))\psi(\q)}{\gcd(\q)}}^m = \infty.
\end{cases}
    \]
\end{thm}

In particular, this theorem contains in it the inhomogeneous Duffin--Schaeffer conjecture for $n$-by-$m$ systems of linear forms with $m\geq3$. For a constant $\y$, Theorem \ref{thm:bootstrapped} follows immediately from~\cite[Theorem~2.5]{ramirez2024metric}, where it is shown that the inhomogeneous Duffin--Schaeffer conjecture in dimension $m$ implies the inhomogeneous Duffin--Schaeffer conjecture for $n$-by-$m$ forms. A full proof of Theorem \ref{thm:bootstrapped} is given in Section~\ref{sec:bootstrap}.

In~\cite{AR_2023} it is shown that Theorem \ref{thm:bootstrapped} holds univariately---meaning the function $\psi$ depends only on $\abs{\q}$---and with a constant $\y$, for all $m\in\NN$ and $n\geq 3$. It is conjectured there that the same should hold for $n\leq 2$. The $n=1$ parts of that conjecture are of course the usual $m$-dimensional inhomogeneous Duffin--Schaeffer conjecture. Theorem \ref{thm:bootstrapped} answers the $m\geq 3$ parts of~\cite[Conjecture~1]{AR_2023} in the even stronger sense that allows multivariate functions---functions $\psi$ that depend on $\q$ and not only $\abs{\q}$.

Notice that the target's movement in Theorem \ref{thm:bootstrapped} is not entirely free. We leave open the question of whether one can relax the theorem so that it applies with an arbitrary $\y = (\y_\q)_{\q\in\ZZ^n}$, as well as that same question in the dimensions where $m\leq 2$. Theorem \ref{thm:weakcounterex} implies that the arbitrary moving target cannot possibly work in the $m=1$ case, even with $n\geq 2$. There, one can simply recreate the construction of Theorem \ref{thm:weakcounterex} by supporting $\psi:\ZZ^n\to[0,\infty)$ on a single ray in $\ZZ^n$.\\

Quite interestingly, it seems that for various setups, the boundary of what existing methods can achieve is the case where $m \geq 3$ (or more generally, $mn \geq 3$ in the linear forms case). The inhomogeneous Khintchine--Groshev Theorem holds without monotonicity for $mn \geq 3$, also allowing a moving target~\cite{allen2021independence}. Similarly, the univariate inhomogeneous Duffin--Schaeffer conjecture for systems of linear forms holds for $n \geq 3$~\cite{AR_2023} and, by Theorem \ref{thm:movingdsc}, the more general multivariate version of it holds for $m \geq 3$. Both cases allow a moving target, albeit restricted movement in the multivariate case.\\

Note that for $m = n= 1$, Khintchine's Theorem fails without the monotonicity assumption, not only in the homogeneous case~\cite{duffinschaeffer}, but also for any fixed inhomogeneous parameter~\cite{ds_counterex}. Among the few positive results concerning the removal of monotonicity in low dimensions in the inhomogeneous setup, there is the nonmonotonic Khintchine--Groshev for non-Liouville $y$ in $(m,n) = (1,2)$~\cite{H23}, as well as for $(m,n) = (1,1)$ assuming restrictions on $\psi$ and $y$ \cite{Yu}. The proofs of all these results make crucial use of the fact that the target is not allowed to move.

\begin{table}[h!]
  \centering
  \begin{center}
    \begin{tabular}{|c||c||c||c|}
      \hline
      \multicolumn{4}{|c|}{\textbf{Khintchine--Groshev-type results and problems}}\\
      \hline\hline
      $n\times m$ & homogeneous & inhomogeneous & moving target \\
      \hline\hline
      \multirow[c]{2}{*}[0in]{\makecell{$n=1$ \\ $m=1$}} & \scriptsize Khintchine~\cite[1924]{K24}  & \scriptsize Szüsz~\cite[1958]{SzuszinhomKT} &  \scriptsize Conjecture~\ref{q:1d} \\ \cline{2-4}
      & \scriptsize \makecell{monotonicity required: \\ Duffin--Schaeffer~\cite[1941]{duffinschaeffer}} & \scriptsize\makecell{monotonicity required: \\ Ram{\'i}rez~\cite[2017]{ds_counterex}} &  \scriptsize\makecell{monotonicity required: \\ $\longleftarrow$ follows~\cite{duffinschaeffer,ds_counterex}} \\
      \hline\hline
       \multirow[c]{2}{*}[0in]{\makecell{$n=1$ \\ $m = 2$}} & \scriptsize Khintchine~\cite[1926]{K26} & \scriptsize Schmidt~\cite[1964]{S64} & \scriptsize Theorem~\ref{thm:movingtarget} \\
      \cline{2-4}
      & \scriptsize\makecell{nonmonotonic $\psi$: \\ Gallagher~\cite[1965]{G65}} & \scriptsize\makecell{nonmonotonic $\psi$: \\ \cite[Conjecture~1]{allen2021independence}} & \scriptsize\makecell{nonmonotonic $\psi$: \\ \cite[Conjecture~2]{allen2021independence}}\\
      \hline\hline
      \multirow[c]{2}{*}[-0.1in]{\makecell{$n=1$ \\ $m\geq 3$}} & \scriptsize Khintchine~\cite[1926]{K26} & \scriptsize Schmidt~\cite[1964]{S64} & \scriptsize Theorem~\ref{thm:movingtarget} \\
      \cline{2-4}
      & \scriptsize\makecell{nonmonotonic $\psi$: \\ Gallagher~\cite[1965]{G65}} & \scriptsize\makecell{nonmonotonic $\psi$: \\ Yu~\cite[2021]{Yu}} & \scriptsize\makecell{nonmonotonic $\psi$: \\ Yu~\cite[2021]{Yu} \\ see also~\cite[Theorem~2]{allen2021independence}}\\
      \hline\hline
      
      \multirow[c]{3}{*}[-0.15in]{\makecell{$n = 2$ \\ $m = 1$}} &  \scriptsize Khintchine--Groshev~\cite[1938]{Groshev}&  \scriptsize Sprind\v{z}uk~\cite[Theorem~15]{Sprindzuk}
                                    & \scriptsize \makecell{Sprind\v{z}uk~\cite[Theorem~15]{Sprindzuk} \\ restricted movement}\\
      \cline{2-4}
      & \scriptsize\makecell{nonmonotonic univariate $\psi$: \\ Beresnevich--Velani~\cite[2010]{BV10}}& \scriptsize \makecell{nonmonotonic univariate $\psi$: \\ 
      \cite[Conjecture~1]{allen2021independence} \\ Hauke~\cite[2024]{H23}, $y$ non-Liouville} & \scriptsize\makecell{nonmonotonic univariate $\psi$: \\ 
      \cite[Conjecture~2]{allen2021independence}}\\
      \cline{2-4}
      &  \multicolumn{3}{c|}{\scriptsize\makecell{multivariate $\psi$: \\ ruled out by~\cite{duffinschaeffer,ds_counterex}}}\\
      \hline\hline
      \multirow[c]{3}{*}[-0.2in]{\makecell{$n=2$ \\ $m\geq 2$}} &  \scriptsize Khintchine--Groshev~\cite[1938]{Groshev}&  \multicolumn{2}{c|}{\scriptsize Sprind\v{z}uk~\cite[Theorem~15]{Sprindzuk}, restricted movement}
                                    \\
      \cline{2-4}
      & \scriptsize\makecell{nonmonotonic univariate $\psi$: \\ Beresnevich--Velani~\cite[2010]{BV10}}& \multicolumn{2}{c|}{\scriptsize \makecell{nonmonotonic univariate $\psi$: Allen--Ram{\'i}rez~\cite{allen2021independence}}} \\
      \cline{2-4}
      &  \scriptsize\makecell{multivariate $\psi$: \\ Beresnevich--Velani~\cite[2010]{BV10}} & \scriptsize\makecell{multivariate $\psi$: \\ Ram{\'i}rez~\cite[2024]{ramirez2024metric}, $m \geq 3$ \\ $m=2$ open} & \scriptsize\makecell{multivariate $\psi$: \\
      open}\\
      \hline\hline
      \multirow[c]{3}{*}[-0.25in]{\makecell{$n\geq 3$ \\ $m\geq 1$}} &  \scriptsize Khintchine--Groshev~\cite[1938]{Groshev}&  \multicolumn{2}{c|}{\scriptsize Sprind\v{z}uk~\cite[Theorem~15]{Sprindzuk}, restricted movement}\\
      \cline{2-4}
      & \scriptsize\makecell{nonmonotonic univariate: \\ Schmidt~\cite{Schmidtcanada}} & \scriptsize\makecell{nonmonotonic univariate $\psi$: \\ Sprind\v{z}uk~\cite[Theorem~15]{Sprindzuk}}& \scriptsize\makecell{nonmonotonic univariate $\psi$: \\ Allen--Ram{\'i}rez~\cite{allen2021independence}}\\
      \cline{2-4}
      &  \scriptsize\makecell{multivariate $\psi$: \\ Beresnevich--Velani~\cite[2010]{BV10}, $m\geq 2$ \\ $m=1$ ruled out} & \scriptsize\makecell{multivariate $\psi$: \\ Ram{\'i}rez~\cite[2024]{ramirez2024metric}, $m \geq 3$ \\ $m=2$ open, $m=1$ ruled out} & \scriptsize\makecell{multivariate $\psi$: \\
      $m=1$ ruled out \\ $m\geq 2$ open}\\
      \hline
    \end{tabular}
  \end{center}
  
  \caption{A table of what is known and what is open for
    Khintchine--Groshev-type problems. For $\psi:\ZZ^n\to[0,\infty)$
    satisfying the appropriate divergence condition, does $W_{n,m}(\y,\psi)$ always have full
    measure? Does it require monotonicity? When $n\geq 2$, is it known univariately, meaning for functions only depending on $\abs{\q}$?}
    \label{tableKG}
\end{table}

\begin{table}[h!]
  \centering
  \begin{center}
    \begin{tabular}{|c||c||c||c|}
      \hline
      \multicolumn{4}{|c|}{\textbf{Duffin--Schaeffer-type results and problems}}\\
      \hline\hline
      $n\times m$ & homogeneous & inhomogeneous & moving target \\
      \hline\hline
      \multirow[c]{3}{*}[-0.2in]{\makecell{$n=1$ \\ $m=1$}} & \multirow[c]{3}{*}[-0.2in]{\scriptsize \makecell{Duffin--Schaeffer conjecture~\cite[1941]{duffinschaeffer} \\  Koukoulopoulos--Maynard~\cite[2020]{KMDS}}} & \scriptsize \makecell{\begin{tabular}{c} Conjectured in~\cite{ds_counterex}, open \\ 
                                                                                                                                                                                                                \end{tabular}}&  \scriptsize\makecell{ruled out \\ by Theorem~\ref{thm:weakcounterex}}\\
      \cline{3-4}
       & & \scriptsize\makecell{weak version:~\cite[Conjecture~1.22]{CT24} \\ Chow--Technau~\cite[2024]{CT24}, cases 
      \\ Beresnevich--Hauke--Velani~\cite[2024]{BHV24}, cases} & \scriptsize \makecell{weak version: \\ Question \ref{quest}} \\
      \cline{3-4}
       & & \scriptsize\makecell{shift-reduced discussions: \\  Schmidt~\cite[1964]{S64}, \\ Chow--Technau~\cite[2024]{CT24}, \\ Beresnevich--Hauke--Velani~\cite[2024]{BHV24}}& \scriptsize open\\
      \hline\hline 
      \multirow[c]{2}{*}[0in]{\makecell{$n=1$ \\ $m = 2$}} & \multirow[c]{3}{*}[-0.1in]{\scriptsize \makecell{Sprind\v{z}uk~\cite[Conjecture]{Sprindzuk} \\ Pollington--Vaughan~\cite[1990]{PV90}}} & \scriptsize\makecell{open} & \scriptsize \makecell{Conjecture \ref{conj:weakdscmoving}}\\
      \cline{3-4}
       & & \scriptsize weak version: open & \scriptsize\makecell{weak version: open}\\
      \cline{1-1}\cline{3-4}
      \makecell{$n=1$ \\ $m\geq 3$} & & \multicolumn{2}{c|}{\scriptsize Theorem~\ref{thm:movingdsc}} \\
      \hline\hline
      \multirow[c]{2}{*}[0in]{\makecell{$n = 2$ \\ $m \leq 2$}} &  \multirow[c]{6}{*}[0in]{\scriptsize\makecell{Beresnevich--Bernik--\\Dodson--Velani~\cite[Conjecture]{BBDV} \\ Ram{\'i}rez~\cite[2024]{ramirez2023duffinschaeffer}}} &  \scriptsize \scriptsize \makecell{see \cite[Conjecture~2.4]{ramirez2024metric}. \\ implied by corresp. $n=1$ case~\cite{ramirez2024metric}}
           & \multirow[c]{2}{*}[0in]{\scriptsize open}\\
      \cline{3-3}
      & &  \scriptsize univariate version: \cite[Conjecture~1]{AR_2023} &  \\
      \cline{1-1}\cline{3-4}
      \multirow[c]{2}{*}[0in]{\makecell{$n = 2$ \\ $m \geq 3$}} &   &  \multirow[c]{2}{*}[0in]{\scriptsize\makecell{Theorem~\ref{thm:bootstrapped}}}
           & \multirow[c]{3}{*}[0in]{\scriptsize\makecell{mostly 
 open, \\ Theorem~\ref{thm:bootstrapped}: \\ restricted movement \\ $m=1,2$ open}}\\
      & &  &  \\
      \cline{1-1}\cline{3-3}
      \multirow[c]{2}{*}[0in]{\makecell{$n\geq 3$ \\ $m\geq 1$}} &   &  \scriptsize\makecell{Theorem~\ref{thm:bootstrapped}, $m\geq 3$ \\ $m=1,2$ open}
           & \\
      \cline{3-4}
      & &  \multicolumn{2}{c|}{\scriptsize \makecell{univariate version:\\
      Allen--Ram{\'i}rez~\cite[2023]{AR_2023}}}   \\
      \hline
    \end{tabular}
  \end{center}
  
  \caption{A table of what is known and what is open for
    Duffin--Schaeffer-type problems. For $\psi:\ZZ^n\to[0,\infty)$
    satisfying the appropriate divergence condition, does $W_{n,m}'(\y,\psi)$ always have full
    measure? And weak versions: with the same divergence condition, does $W_{n,m}(\y,\psi)$ always have full measure?}
    \label{tableDS}
\end{table}

\subsubsection*{Khintchine's Theorem with a moving target}

Returning attention to Khintchine's Theorem (with the monotonicity assumption), the question arises whether a moving target version holds. It turns out that when $m\geq 2$, the inhomogeneous version of Khintchine's Theorem can be stated with a moving target. 

\begin{thm}\label{thm:movingtarget}
  Let $m \geq 2$. For any sequence
  $\y = (\y_q)_{q\in \N}\subset \RR^m$ and monotonic function $\psi:\NN\to[0,\infty)$ such that $\sum_{q \in \N}\psi(q)^m$
  diverges, the set $W_m(\y,\psi)$
  has full Lebesgue measure in $[0,1]^m$. 
\end{thm}

The cases where $m \geq 3$ are implicit in the work of Yu~\cite{Yu} and made explicit in~\cite[Theorem~2]{allen2021independence}. We can find no proof in the literature for the $m=2$ case
of Theorem \ref{thm:movingtarget}, so a short proof is provided in Appendix~\ref{sec:two},
for completeness. We conjecture that Khintchine's original one-dimensional theorem is also true with a moving target
$y= (y_q)_{q\in\NN}$.

\begin{conj}\label{q:1d}
Theorem \ref{thm:movingtarget} holds for $m=1$.
\end{conj}

Conjecture \ref{q:1d} would follow immediately from a positive answer to Question \ref{conj:weakdscmoving} (the weak Duffin--Schaeffer conjecture with a moving target), but it might be the case that the monotonicity condition assumed in Conjecture \ref{q:1d} makes a direct proof more attainable, or possible, while Question \ref{conj:weakdscmoving} might turn out too much to ask for.

Regarding monotonicity assumptions, it is found in~\cite{allen2021independence,Yu} that Theorem \ref{thm:movingtarget} does not need a monotonic $\psi$ in dimensions $m\geq3$. The proof we give for the $m=2$ case makes use of the monotonicity of $\psi$. In~\cite[Conjecture~2]{allen2021independence} it is conjectured that the $m=2$ case of Theorem \ref{thm:movingtarget} also holds without monotonicity.

\subsubsection*{Hausdorff measure statements}

In 2006, Beresnevich and Velani proved the mass transference principle~\cite{BVmassTP}, a result that allows one to deduce Hausdorff measure statements from Lebesgue measure statements. The mass transference principle has enjoyed wide application throughout the years. In its first application, Beresnevich--Velani showed that the $m$-dimensional Duffin--Schaeffer conjecture implies the $m$-dimensional Duffin--Schaeffer conjecture for Hausdorff measures, for all $m\geq 1$. The exact same argument works in the inhomogeneous setting, leading immediately to the following.

\begin{thm}
    Let $m\geq 3$, and $f:\RR_+\to \RR_+$ a dimension function so that $x^{-m} f(x)$ is monotonic. Then for arbitrary sequences $\y = (\y_q)_{q\in \N}$ and functions $\psi:\NN\to[0,\infty)$, we have 
    \begin{equation*}
        \calH^f\left(W_{n,m}'(\y,\psi)\right) =
        \begin{cases}
            0 &\textrm{if } \sum\limits_{q\in \N} \varphi(q)^m f\parens*{\frac{\psi(q)}{q}} <\infty, \\
            \calH^f\left([0,1]^m\right) &\textrm{if } \sum\limits_{q\in \N} \varphi(q)^m f\parens*{\frac{\psi(q)}{q}} =\infty,
        \end{cases}
    \end{equation*}
    where $\calH^f$ denotes the Hausdorff $f$-measure on $\RR^m$.
\end{thm}

\section{Proof of Theorem \ref{thm:movingdsc}}

\subsection{Proof strategy and insights}

The statement of Theorem \ref{thm:movingdsc} can be rephrased in measure-theoretic language by asking about the Lebesgue measure of $\limsup_{q \to \infty} \Aqm$ where
$\Aqm$ denotes the event
\begin{equation*}
\lvert q\mathbf{x} - \ba - \byq \rvert \leq \psi(q),\quad \textrm{ for some $\ba$ with }\quad \; (a_1,q) = \ldots = (a_m,q)= 1.
\end{equation*}
As is usual in metric Diophantine approximation, the argument is carried out by bounding the pairwise dependence of the events $\Aqm \cap \Arm$ on average from above, in order to apply some refined form of the divergence Borel--Cantelli Lemma (in our case Proposition \ref{BC5}). Since the events under consideration here are indicators of unions of intervals, these estimates are often called overlap estimates. Note that the indicators are supported on small (hyper-)cubes, which reduces the task to obtaining overlap estimates in dimension $m = 1$. Therefore, the overlap estimate established here can also be seen as a step towards the ($1$-dimensional) inhomogeneous Duffin--Schaeffer conjecture. However, the error terms we obtain are too big to obtain anything meaningful for the $1$-dimensional setting. We remark that our proof does not make use of the recent Koukoulopolus--Maynard method; it is sufficient to apply the ideas established by Pollington--Vaughan.\\

For the overlap estimate itself, we follow the general strategy of Pollington and Vaughan in~\cite[Section 3]{PV90}, borrowing also some ideas of~\cite[Lemma 5]{ABH23} and combine this with new sieve-theoretic input. In comparison to the homogeneous theory established in these preceding articles, an additional challenge arises: The inhomogeneous target leads to sieve estimates on short intervals of length $D$, where $D$ varies for different pairs of denominators $(q,r)$, depending on their arithmetic structures as well as their corresponding $\psi$-weights. In the extreme case where $D < \tfrac{1}{2}$, the coprimality assumption no longer prevents overlaps (as it does in the homogeneous setup), and this yields an error term that seems unavoidable without changing the sets. If the sifting range $D$ is long enough, we can overcome this by some adapted sieve estimates at the cost of an additional error term, and this error can be shown to be dominated by the (seemingly unavoidable) error term arising from the case of $D < \tfrac{1}{2}$. Fortunately, it turns out that in dimension $m \geq 3$, the extra errors can be handled by weighted gcd-sums. Note that Theorem \ref{thm:weakcounterex} shows that in dimension $m=1$, not only the classical $L^2$-approach fails, but indeed, there is no possibility of circumventing this issue for a moving target. Thus, as already mentioned, if there is a possibility to prove the inhomogeneous Duffin--Schaeffer conjecture, the proof cannot build solely on controlling the overlaps of pairs $(q,r)$ with arbitrary target; rather, one needs to take into account the special shape given by the fixed inhomogeneous parameter.

\subsection{The overlap estimate} \label{sec_3}

\begin{lem}[Overlap estimate in one dimension] \label{lemma_over}
	Let $\psi: \N \to [0,\infty)$, $q \neq r$ be positive integers, $y,z \in \mathbb{R}$ and let
 \[\begin{split}\Aq &:= \bigcup_{(a,q) = 1} \left[\frac{a + y}{q} - \frac{\psi(q)}{q}, \frac{a + y}{q} + \frac{\psi(q)}{q}\right] \cap [0,1), \\
 \Ar &:= \bigcup_{(b,r) = 1} \left[\frac{a + z}{r} - \frac{\psi(r)}{r}, \frac{a + z}{r} + \frac{\psi(r)}{r}\right] \cap [0,1).\end{split}
 \]
Writing $D(q,r) := 2\lcm(q,r)\max\left\{\frac{\psi(q)}{q},\frac{\psi(r)}{r}\right\}$, we have
	\begin{equation} \label{overlapclaim}
	\overlap \ll \mathds{1}_{[D(q,r) > 1]}\frac{\psi(q)\varphi(q)}{q}\frac{\psi(r)\varphi(r)}{r}\prod_{\substack{p \mid \frac{qr}{(q,r)^2}\\ p > D(q,r)}}\left(1 + \frac{1}{p}\right) + \varphi(\gcd(q,r))\min\left\{\frac{\psi(q)}{q},\frac{\psi(r)}{r}\right\},
	\end{equation}
	with an absolute implied constant.
\end{lem}

The proof of Lemma \ref{lemma_over} contains the main new ingredients needed in order to prove Theorem \ref{thm:movingdsc}. To establish the statement, we apply results from sieve theory, which we will do in the form of the fundamental lemma of sieve theory, using the formulation of~\cite[Theorem 18.11]{K19}.
\begin{lem}[Fundamental lemma of sieve theory] \label{lemma_fund}
	Let $(a_n)_{n \in \N}$ be non-negative reals with $\sum_{n=1}^\infty a_n < \infty$. Let $\mathcal{P}$ be a finite set of primes, and write $P = \prod_{p \in \mathcal{P}} p, y = \max \mathcal{P}$, and $E_d = \sum_{n \equiv 0 \mod d} a_n$. Assume that there exists a multiplicative function $g$ such that $g(p) < p$ for all $p \in \mathcal{P}$, a real number $x$, and positive constants $\kappa,C$ such that
	$$
	E_d =: x \frac{g(d)}{d} + r_d, \qquad d \mid P,
	$$
	and
	$$
	\prod_{p \in (y_1, y_2] \cap \mathcal{P}} \left( 1 - \frac{g(p)}{p} \right)^{-1} < \left( \frac{\log y_2}{\log y_1} \right)^\kappa \left(1 + \frac{C}{\log y_1} \right), \qquad 3/2 \leq y_1 \leq y_2 \leq y.
	$$
	Then, uniformly in $u \geq 1$ we have
	$$
	\sum_{(n,P)=1} a_n = \left( 1 + O ( u^{-u/2} ) \right) x \prod_{p \in \mathcal{P}} \left(1  -\frac{g(p)}{p} \right) + O \Bigg( \sum_{\substack{d \leq y^u,\\d \mid P}} |r_d| \Bigg).
	$$
\end{lem}

\begin{proof}[Proof of Lemma \ref{lemma_over}]

Following ~\cite[Section 3]{PV90}, we set 
	\begin{align*}
		\delta := 2\min \left( \frac{\psi(q)}{q}, \frac{\psi(r)}{r} \right), \qquad  \qquad \Delta := 2\max \left( \frac{\psi(q)}{q}, \frac{\psi(r)}{r} \right),\\
	\end{align*}
	and bound the overlap straightforwardly by
	\[
 \begin{split}
	\overlap &\ll \delta \underset{\left| \frac{a}{q} - \frac{y_q}{q} - \left(\frac{b}{r} - \frac{y_r}{r}\right)\right| \leq \Delta
    }{\sum_{\substack{1 \leq a \leq q,\\ (a,q) = 1}}  \sum_{\substack{1 \leq b \leq r,\\ (b,r) = 1}}} 1= \delta \sum_{c \in \mathbb{Z} \cap [X,Y]}f(c)
    \end{split}
	\]
where 
\[[X,Y] = \left[X(q,r,y_q,y_r), Y(q,r,y_q,y_r)\right] := \left[-\tfrac{D}{2},\tfrac{D}{2}\right] + \frac{r}{(q,r)}y_q - \frac{q}{(q,r)}y_r,\]
and
\[f(c) = f_{q,r}(c) := \#\left\{(a,b) \in \mathbb{Z}_q^{*} \times \mathbb{Z}_r^{*}: \frac{a}{q} - \frac{b}{r} = \frac{c}{\lcm(q,r)}\right\}.\]
	For any prime $p$, let $u = u(p,q)$ and $v = v(p,r)$ be defined by $q = \prod_p p^u$ and $r = \prod_p p^v$, and let
 $$ \ell = \prod_{p:~u=v} p^u, \qquad m = \prod_{p:~u \neq v} p^{\min(u,v)}, \qquad n = \prod_{p:~u \neq v} p^{\max(u,v)}.$$
	Following the argument of~\cite[p.195]{PV90} (an application of the Chinese remainder theorem, together with an elementary counting argument) leads to

 \begin{equation}\label{coprime_count}f(c) =  \mathds{1}_{[(c,n) = 1]}\varphi(m) \ell \prod_{p \mid (\ell,c)} \left(1 - \frac{1}{p} \right) \prod_{\substack{p \mid \ell, \\p \nmid c}} \left( 1 - \frac{2}{p} \right) \ll \mathds{1}_{[(c,n) = 1]}
 \varphi(m) \frac{\varphi(\ell)^2}{\ell} \prod_{p \mid (\ell,c)} \left(1 + \frac{1}{p} \right),
	\end{equation}
 where $\mathds{1}$ denotes the indicator function.

If $D = Y-X \leq 100$, then we use the pointwise upper bound
$f(c) \leq \varphi(m)\varphi(\ell) = \varphi(\gcd(q,r))$ to obtain 
\[\overlap \ll \delta (D+1) \varphi(\gcd(q,r)) \ll \delta\varphi(\gcd(q,r)),\]
which proves \eqref{overlapclaim} in case of $D \leq 100$. Thus, from now on we can assume that $D = Y-X> 100$.\\

In this case, we note that
 \begin{align}\overlap &\ll \delta \varphi(m) \frac{\varphi(\ell)^2}{\ell} 
 \sum_{\substack{c \in \mathbb{Z} \cap [X,Y]\\(c,n) = 1}} \prod_{p \mid (\ell,c)} \left(1 + \frac{1}{p} \right)
 \leq \delta \varphi(m) \frac{\varphi(\ell)^2}{\ell} 
\sum_{\substack{c \in \mathbb{Z} \cap [X,Y]\\(c,n) = 1}}  \sum_{d \mid (\ell,c)}\frac{1}{d}.\label{before_sieve}
 \end{align}
Next, we fix a parameter $T$ determined later. Using $(\ell,n) = 1$, we obtain

\[\begin{split}\sum_{\substack{c \in \mathbb{Z} \cap [X,Y]\\(c,n) = 1}}   \sum_{d \mid (\ell,c)}\frac{1}{d}
&
= \sum_{\substack{d \mid \ell\\ d \leq T}} \frac{1}{d}\sum_{\substack{c \in \mathbb{Z} \cap \left[\tfrac{X}{d},\tfrac{Y}{d}\right]\\(c,n) = 1}} 1 +  \sum_{\substack{d \mid \ell\\ d > T}} \frac{1}{d}\sum_{\substack{c \in \mathbb{Z} \cap \left[\tfrac{X}{d},\tfrac{Y}{d}\right]\\(c,n) = 1}} 1 =: S_1 + S_2.
\end{split}
\]

For $S_1$, we perform for every $d \leq T$ fixed an upper-bound sieve. Writing 
$n^*$ for the $T$-smooth part of $n$, we observe that for any $e \mid n^*$ we get
	\begin{eqnarray*}
		\sum_{\substack{c \in \mathbb{Z} \cap \left[\tfrac{X}{d},\tfrac{Y}{d}\right] \\ c \equiv 0 \mod e}} 1 & = \frac{Y-X}{de} + O(1).
	\end{eqnarray*}
	By an application of Lemma~\ref{lemma_fund} (with $\mathcal{P}$ the set of prime divisors of $n^*$, $\max \mathcal{P} \le T$ and $|r_e| \ll 1, u = 1$), this shows

 \[
\sum_{\substack{c \in \mathbb{Z} \cap \left[\tfrac{X}{d},\tfrac{Y}{d}\right]\\(c,n) = 1}} 1
\ll \frac{D}{d} \prod_{p \mid n^*}\left(1 - \frac{1}{p}\right) + O \Bigg( \sum_{\substack{e \leq T,\\e \mid P}} 1 \Bigg) = \frac{D}{d} \frac{\varphi(n^*)}{n^*} + O(T),
 \]
 which after summing over $d$ implies

 \[
 S_1 \ll D\frac{\varphi(n^*)}{n^*}\Bigg(\sum_{\substack{d \mid \ell\\ d \leq T}} \frac{1}{d^2}\Bigg) + O(T \log T).
 \]
For $S_2$, we simply drop the condition $(c,n) = 1$ and deduce 
\[\begin{split}S_2 \leq \sum_{\substack{d \mid \ell\\ d > T}} \frac{1}{d}\sum_{\substack{c \in \mathbb{Z} \cap \left[\tfrac{X}{d},\tfrac{Y}{d}\right]}} 1 
\ll \sum_{\substack{d \mid \ell\\ d > T}} \left(\frac{D}{d^2} + \frac{1}{d}\right)
\ll \frac{D}{T} + \sum_{\substack{d \mid \ell}} \frac{1}{d}.
\end{split}
\]
Using $\prod_{p \leq T}\left(1 - \frac{1}{p}\right) \gg \frac{1}{\log T}$ and $\sum_{\substack{d \mid \ell}} \frac{1}{d} \ll \frac{\ell}{\varphi(\ell)}$, we obtain

\[S_1 + S_2 \ll 
D\frac{\varphi(n^*)}{n^*}\left(1 + O\left(\frac{T (\log T)^2}{D} + \frac{\log T}{T}\right)\right) + 
\frac{\ell}{\varphi(\ell)}.
\]
Choosing $T = D^{1/3}$ finally shows that

\begin{equation}\label{S1S2}
S_1 + S_2 \ll D\frac{\varphi(n)}{n}\prod_{\substack{p \mid n \\ p > D^{1/3}}}\left(1 + \frac{1}{p}\right) + \frac{\ell}{\varphi(\ell)}
\ll D\frac{\varphi(n)}{n}\prod_{\substack{p \mid n \\ p > D}}\left(1 + \frac{1}{p}\right) + \frac{\ell}{\varphi(\ell)},
\end{equation}
where we used Mertens' Theorem in the form of 
$\prod_{\substack{D^{1/3}<  p < D}}\left(1 + \frac{1}{p}\right) \ll 1.$
Plugging \eqref{S1S2} into \eqref{before_sieve}, using
\[\varphi(m)\varphi(\ell) = \varphi(\gcd(q,r)), \quad \varphi(m)\varphi(\ell)^2\varphi(n) = \varphi(q)\varphi(r), \quad m \ell D\delta = \psi(q)\psi(r),\]
we get 
\[\begin{split}\overlap
&\ll \delta \varphi(m) \frac{\varphi(\ell)^2}{\ell} \Bigg(D\frac{\varphi(n)}{n}\prod_{\substack{p \mid n \\ p > D}}\Big(1 + \frac{1}{p}\Big) + \frac{\ell}{\varphi(\ell)}\Bigg)
\\&= \frac{\psi(q)\varphi(q)}{q}\frac{\psi(r)\varphi(r)}{r}\prod_{\substack{p \mid n \\ p > D}}\left(1 + \frac{1}{p}\right) + \delta\varphi(m)\varphi(\ell)
\\&= \frac{\psi(q)\varphi(q)}{q}\frac{\psi(r)\varphi(r)}{r} \prod_{\substack{p \mid \frac{qr}{(q,r)^2} \\ p > D}}\left(1 + \frac{1}{p}\right) + \varphi(\gcd(q,r))\min\left\{\frac{\psi(q)}{q},\frac{\psi(r)}{r}\right\},
\end{split}
\]
which concludes the proof.
\end{proof}

\subsection{Finishing the proof}

In inhomogeneous Diophantine approximation, there are no $0$-$1$-laws analogous to the homogeneous $0$-$1$ laws of Cassels and Gallagher~\cite{C01,G01}. The absence of $0$-$1$-laws for the inhomogeneous setup is usually overcome by working locally. One uses the regularity properties of $A_q^y$ to establish quasi-independence restricted to an arbitrary small set $U$, and then an argument based on the Lebesgue density theorem allows one to conclude full measure. This was carried out in various related articles, e.g. ~\cite{allen2021independence,AR_2023,CT24,ramirez2023duffinschaeffer,ramirez2024metric} and can also be applied here. However, to streamline the argument, we will instead apply a recent refinement of the Borel--Cantelli Lemma~\cite{BHV24} that avoids the use of Lebesgue density arguments.

\begin{prop}[{~\cite[Theorem 5]{BHV24}}]\label{BC5}
Let $\mu$ be a doubling Borel regular probability measure on a metric space $X$. Let $(E_i)_{i\in\N}$ be a sequence of $\mu$-measurable subsets of $X$. Suppose that 

\begin{equation}\label{eqn01}
\sum_{i=1}^\infty \mu(E_i)=\infty
\end{equation}
and that there exists a constant $C>0$ such that
\begin{equation}\label{eqn02}
\sum_{s,t=1}^Q  \mu(E_s\cap E_t)\le C\left(\sum_{s=1}^Q  \mu(E_s)\right)^2\quad\text{for infinitely many $Q\in\N$\,.}
\end{equation}

 In addition, suppose that for any $\delta>0$ and any closed ball $B$ centred at $\supp\mu$ there exists $i_0 =i_0 (\delta, B) $ such that for all $i  \geq i_0$  for all large enough $i$
\begin{equation}\label{vb89}
\mu\left(B\cap E_i\right)\le (1+\delta)\mu\left(B\right)\mu(E_i)\,.
\end{equation}
Then $\mu(E_\infty)=1$.
\end{prop}

We would like to replace $\lambda_m(\Aqm)$ by $\Big(2\frac{\psi(q)\varphi(q)}{q}\Big)^m$ in order to be able to satisfy \eqref{eqn01} and \eqref{eqn02} for our setup.
This holds trivially if $\psi(q) \leq \frac{1}{2}$, but for larger values of $\psi$, this is in general no longer true due to possible overlaps of the intervals contributing to $\Aqm$, which complicates matters. Similarly, \eqref{vb89} follows for $\psi(q) \leq \frac{1}{2}$ by inheriting the equidistribution property from $\mathbb{Z}_q^{*}$, but for large values of $\psi$, some information about the gap sizes of $\mathbb{Z}_q^{*}$ are needed.

However, this issue can be treated with the following statement, provided $\psi$ is not too large. This was recently proven in~\cite{BHV24} in a more general setup, with the special case considered here refining the method established in~\cite[Chapter 2]{PV90}. 

\begin{prop}[special case of {\cite[Proposition 2]{BHV24}}]\label{prop:mixing}
    Let $\psi: \N \to [0,\infty)$ and suppose that
\begin{equation}\label{decay_large_psi}\lim_{\substack{q \to \infty\\  {\psi(q)\ge1/2}}} \frac{\varphi(q)\psi(q)}{q} = 0.\end{equation}
    Let $y = (y_q)_{q \in \mathbb{N}} \subset \R$ be arbitrary and $\Aq$ as in Lemma \ref{lemma_over}. 
     Then for any $0 \leq x < y \leq 1$,
\begin{equation}\label{ud_large_psi}
 \lim_{\substack{q \to \infty\\ {\psi(q) \neq 0}}} \frac{\lambda_1\left(\Aq \cap [x,y] \right)}{\lambda_1(\Aq)} =\lambda_1([x,y]).
\end{equation}
Furthermore, we have
\begin{equation}\label{meas_aq}\lambda_1(\Aq) \gg \frac{\psi(q)\varphi(q)}{q},\end{equation}
with an absolute implied constant.
\end{prop}

Note that both Proposition \ref{prop:mixing} and Lemma \ref{lemma_over} only consider the $1$-dimensional setup. The following is a simple corollary that generalizes the above to the multidimensional variants.
\begin{cor}\label{multidim_statements}
Let $m \geq 1$, $\by = (\byq)_{q \in \N} \subset \R^m$ arbitrary and $\psi: \N \to [0,\infty)$ such that \eqref{decay_large_psi} holds. Then for $q \neq r$, we have

\begin{equation}
    \label{multidim_overlap}
	\lambda_m( \Aqm \cap \Arm) 
 \ll_m M(q,r)^m + \left(\varphi(\gcd(q,r))\frac{\psi(q)}{q}\right)^m,
\end{equation}
where
\[M(q,r) = M(\psi,q,r) := \mathds{1}_{[D(q,r) \geq 1]}\frac{\psi(q)\varphi(q)}{q}\frac{\psi(r)\varphi(r)}{r} \prod_{\substack{p \mid \frac{qr}{(q,r)^2}\\ p > D(q,r)}}\left(1 + \frac{1}{p}\right),\]
with $D$ being as in Lemma \ref{lemma_over}.
Furthermore, we have
\begin{equation}\label{meas_lower}
    \lambda_m\big(\Aqm \big) \gg_m \measqm,
\end{equation}
and for any axis-parallel rectangle $\mathbf{R} \subset [0,1)^m$,
\begin{equation}\label{multidim_equidistr}\lim_{\substack{q \to \infty\\ {\psi(q) \neq 0}}} \frac{\lambda_m\left(\Aqm \cap \mathbf{R} \right)}{\lambda_m(\Aqm)} =\lambda_m(\mathbf{R}).
\end{equation}
\end{cor}

\begin{proof}
    Writing $\by = (y_1,\ldots,y_m)$, we observe that 
    \begin{equation}
        \label{decomp}
    \Aqm = A_q^{y_1} \times A_q^{y_2} \times \ldots \times A_q^{y_m}
    \end{equation}
    which implies
    \[\lambda_m( \Aqm \cap \Arm) 
    = \lambda_1(A_q^{y_1} \cap A_r^{y_1}) \times \ldots \times \lambda_1(A_q^{y_m} \cap A_r^{y_m}).
    \]
    Thus an application of Lemma \ref{lemma_over} and inequality $(a+b)^m \leq 2^m(a^m + b^m)$ proves 
    \eqref{multidim_overlap}. Using \eqref{decomp} again, \eqref{meas_lower} and \eqref{multidim_equidistr} follow immediately from \eqref{meas_aq} resp \eqref{ud_large_psi}.
\end{proof}

The following lemma was crucial in establishing the (homogeneous) Duffin--Schaeffer conjecture, both in dimension $1$ as well as in dimension $2$, eventually bounding the main term of the overlap estimates on average by the corresponding products.

\begin{lem}\label{pv_estimate}
Let $\psi: \N \to \left[0,\tfrac{1}{2}\right]$
with $\sum_{q \in \mathbb{N}} \left(\frac{\psi(q)\varphi(q)}{q}\right)^m = \infty$ for some
$ m \geq 1$. 
Then for any  $Q$ large enough and $M(q,r)$ as in Corollary \ref{multidim_statements}, we have

    \[
\sum_{q\neq r \leq Q} M(q,r)^m \ll_m \left(\sum_{q \leq Q} \left(\frac{\varphi(q)\psi(q)}{q}\right)^m\right)^2,
    \]
    with the implied constant only depending on $m$.
\end{lem}

\begin{rem}
    Although we only make use of Lemma \ref{pv_estimate} in the case of $m \geq 3$, we decided to state it in the strongest possible version, holding for arbitrary $m \geq 1$. The case where $m = 1$ is implicitly the key ingredient in the work of Koukoulopoulos--Maynard~\cite{KMDS} and can be explicitly found (with even refined bounds) in the proofs of~\cite[Theorem 2]{ABH23}, ~\cite[Corollary 1.3]{HSW} and ~\cite[Theorem 2]{KMY24}. The case where $m \geq 2$ was implicitly proven by Pollington--Vaughan~\cite[Section 4]{PV90}. 
    Although the case where $m \geq 2$ can also be derived by using the quite involved GCD graph approach of Koukoulopoulos--Maynard (for a shorter argument see \cite{HSW}), we include the simpler proof in dimension $m \geq 2$ of Pollington--Vaughan  in a detailed form in Appendix~\ref{sec:PVredo}.
\end{rem}

Finally, we provide a lemma that bounds the second term of our overlap estimates in a sufficient way. We remark that here, $m \geq 3$ is a crucial assumption and the term differs for $m \leq 2$.

\begin{lem}\label{lem:phigcd}
  For $m \geq 3$, we have 
  \begin{equation}
    \label{eq:phigcd}
    \sum_{r = 1}^q \varphi(\gcd(q,r))^m \ll \varphi(q)^m. 
  \end{equation}
  For $m =2$, the sum is bounded by $q^2$. 
\end{lem}

\begin{proof}
  Note that
  \begin{align}
    \sum_{r = 1}^q \varphi(\gcd(q,r))^m &= \sum_{d\mid q}\varphi(d)^m\varphi(q/d) \nonumber \\
    &\leq \varphi(q)\sum_{d\mid q}\varphi(d)^{m-1}.\label{eq:phim-1}
  \end{align}
  For a prime $p$, let $\nu_p$ denote the $p$-adic valuation. Note that
  \begin{align*}
    \sum_{d\mid q}\varphi(d)^{m-1}
    &= \prod_{p \mid q}\sum_{k=0}^{\nu_p(d)}\varphi(p^k)^{m-1} \\
    &= \prod_{p \mid q}\parens*{1 + \sum_{k=1}^{\nu_p(d)}\parens*{p^{k-1}(p- 1)}^{m-1} } \\
    &= \prod_{p \mid q}\parens*{p^{\nu_p(d)-1}(p- 1)}^{m-1} \big(1 + O\big(p^{-(m-1)}\big)\big) \\
    &= \varphi(q)^{m-1} \prod_{p \mid q} \left(1 + O\left(p^{-(m-1)}\right)\right) \\
    &\overset{m \geq 3}{\ll} \varphi(q)^{m-1},
  \end{align*}
which proves the lemma for $m\geq 3$. In the case $m=2$, \eqref{eq:phim-1} becomes 
\begin{equation*}
    \varphi(q)\sum_{d\mid q}\varphi(d) = q\varphi(q) \leq q^2,
\end{equation*}
concluding the proof.
\end{proof}

\begin{proof}[Proof of Theorem \ref{thm:movingdsc}]
  The convergence case is a standard application of the convergence
  Borel--Cantelli lemma. For every $q\in\NN$ we have
  \begin{equation}\label{eq:measAq}
    \lambda_m(\Aqm) \leq \parens*{2\frac{\varphi(q)\psi(q)}{q}}^m,
  \end{equation}
  and so the convergence of $\sum_{q \in \N} \left(\frac{\varphi(q)\psi(q)}{q}\right)^m$ implies that the
  sets $(\Aqm)_{q \in \N}$ have a converging measure sum. Hence, their limsup set
  has measure $0$.\\

  For the divergence case, we claim that it is enough to consider functions $\psi$ satisfying \eqref{decay_large_psi}.
Indeed, assume for the moment that \eqref{decay_large_psi} does not hold. Then
 on letting 
$\psi_1(q) := \min\left\{\psi(q),\frac{q}{\varphi(q)}\right\}$ we  observe that 
$\sum_{q \in \N} \left(\frac{\psi_1(q)\varphi(q)}{q}\right)^m = \infty\,.$
Furthermore, since we assumed that \eqref{decay_large_psi} does not hold, we have that $\limsup_{q \to \infty} \frac{\psi_1(q)\varphi(q)}{q} = \delta$ for some $0 < \delta \leq 1$. Therefore, there exists an increasing sequence of integers $(q_n)_{n \in \N}$ such that for all $n \in \N$ we have
\begin{equation}\label{almost_delta}\frac{\delta}{2} \leq \frac{\psi_1(q_n)\varphi(q_n)}{q_n} \leq 2\delta\,.\end{equation}
Setting
\[\psi_2(q) = \begin{cases} \frac{\psi_1(q_n)}{n^{1/m}} &\text{ if } q = q_n \text { for some $n \in \N$},\\
0 &\text{ otherwise,}
\end{cases}\]
we see that by \eqref{almost_delta} that
\[\sum_{q \in \N} \left(\frac{\psi_2(q)\varphi(q)}{q}\right)^m = \sum_{n \in \N} \left(\frac{\psi_1(q_n)\varphi(q_n)}{q_n}\right)^m \geq \frac{\delta^m}{2^m}\sum_{n \in \N}\frac{1}{n} = \infty\,.
\]
Since $\psi\ge\psi_1\ge\psi_2$ we also have that \[\lambda_m\left(\limsup_{q \to \infty} \Aqm(\psi)\right) \geq \lambda_m\left(\limsup_{q \to \infty} \Aqm(\psi_2)\right),\] and therefore in the proof we could replace $\psi$ with $\psi_2$, and $\psi_2$ clearly satisfies \eqref{decay_large_psi}. This justifies the claim about $\psi$ we made above, and we can assume that $\psi$ satisfies \eqref{decay_large_psi}.\\
  
  Now let
  \begin{equation*}
    I = \set*{q\in\NN : \frac{1}{2} < \psi(q) \leq \frac{q}{2\varphi(q)}}. 
  \end{equation*}
  The rest of this proof is split into two cases.



  \subsection*{Case 1:} The sum $\sum_{q\in I}(\varphi(q)\psi(q)/q)^m$
  diverges. In this case, we may replace $\psi$ with $\psi\mathds{1}_{I}$,
  where $\mathds{1}_{I}$ is the indicator function of $I$. By Corollary \ref{multidim_statements}, we have
  \begin{equation*}
    \lambda_m(\Aqm \cap \Arm)
    \ll \productm\prod_{\substack{p \mid \frac{qr}{(q,r)^2}\\ p > D(q,r)}}\left(1 + \frac{1}{p}\right) + \frac{\psi(q)^m}{q^m}\varphi(\gcd(q,r))^m
    \end{equation*}
  for all $r,q \in I$ with $r < q$. 
Since $\psi(q),\psi(r) > 1$, we have $D(q,r) \geq \sqrt{\frac{qr}{(q,r)^2}}$. Using Mertens' Theorem, we have 
    \[\prod_{\substack{p \mid \frac{qr}{(q,r)^2}\\ p > D(q,r)}}\left(1 + \frac{1}{p}\right)
    \ll \prod_{\substack{D(q,r) < p < D(q,r)^2}}\left(1 + \frac{1}{p}\right) \ll 1.
    \]
 This implies that for $Q>1$,
  \begin{align*}
    \sum_{q,r = 1}^Q\lambda_m(\Aqm \cap \Arm)
    &\overset{\textrm{\eqref{meas_lower}}}{\ll} \sum_{q=1}^Q \lambda_m(\Aqm) + \sum_{q,r=1}^Q\lambda_m(\Aqm)\lambda_m(\Arm) + \sum_{q=1}^Q\frac{\psi(q)^m}{q^m} \sum_{r=1}^q\varphi(\gcd(q,r))^m \\
    &\overset{\textrm{Lem.~\ref{lem:phigcd}}}{\ll} \parens*{\sum_{q=1}^Q\lambda_m(\Aqm)}^2 + \sum_{q=1}^Q\parens*{\frac{\varphi(q)\psi(q)}{q}}^m \\
    &\ll  \parens*{\sum_{q=1}^Q\lambda_m(\Aqm)}^2.
  \end{align*}
 The theorem is proved by applying Proposition \ref{BC5} and \eqref{multidim_equidistr}.

  \subsection*{Case 2:} The sum $\sum_{q\in I}(\varphi(q)\psi(q)/q)^m$
  converges. In this case, we may assume without loss of generality
  that $I=\emptyset$, in other words,
  \begin{equation*}
    \psi(q) \leq \frac{1}{2} \qquad (q\in\NN). 
  \end{equation*}
    Now, for $Q>1$, using Corollary \ref{multidim_statements} and Lemma \ref{pv_estimate}, we have
  \begin{align*}
    \sum_{q,r = 1}^Q\lambda_m(\Aqm \cap \Arm)
    &\ll \sum_{q=1}^Q \lambda_m(\Aqm) + \sum_{\substack{q,r=1 \\ q \neq r}}^Q M(q,r)^m + \sum_{q=1}^Q\frac{\psi(q)^m}{q^m} \sum_{r=1}^q\varphi(\gcd(q,r))^m \\
    &\ll \sum_{q=1}^Q \lambda_m(\Aqm) + \parens*{\sum_{q=1}^Q\parens*{\frac{\varphi(q)\psi(q)}{q}}^m}^2 + \sum_{q=1}^Q\parens*{\frac{\varphi(q)\psi(q)}{q}}^m \\
    &\overset{\textrm{\eqref{meas_lower}}}{\ll} \parens*{\sum_{q=1}^Q\lambda_m(\Aqm)}^2.
  \end{align*}
  Again, the theorem follows by applying Proposition \ref{BC5} and \eqref{multidim_equidistr}.
\end{proof}

\section{Proof of Theorem \ref{thm:weakcounterex}}

The proof relies on some simple lemmas. 

\begin{lem}\label{lem:sumset}
  Let $Q_q'$ denote the reduced fractions with denominator $q$ in
  $\RR/\ZZ$. If $q$ is square-free and $r\mid q$, then
  $Q_r' + Q_{q/r}'= Q_q'$.
\end{lem}

\begin{proof}
  Let $a/r \in Q_r'$ and $b/(q/r) \in Q_{q/r}'$. Then
  \begin{equation*}
    \frac{a}{r} + \frac{b}{q/r} =\frac{c}{q},
  \end{equation*}
  where $c= a(q/r) + rb$. We will show that $c/q \in Q_q'$.

  Suppose $p \mid q$, with $p$ prime. Then, since $q$ is square-free,
  $p$ divides exactly one of $(q/r)$ and $r$.
  \begin{itemize}
  \item If $p$ divides $r$, then it cannot divide $a$, since
    $a/r\in Q_r'$, hence $p$ does not divide $a(q/r)$.
  \item If $p$ divides $(q/r)$, then it cannot divide $b$, since
    $b/(q/r) \in Q_{q/r}'$, hence $p$ does not divide $br$.
  \end{itemize}
  In either case, we have $p\nmid c$. This proves that
  $\gcd(c,q)=1$, hence $c/q \in Q_q'$ as claimed.

  In particular, for any $t\in Q_{q/r}'$, the set $Q_r' + t$ is a
  rotated copy of $Q_r'$ contained in $Q_q'$. Furthermore, for
  $s\in Q_{q/r}'$ with $s\neq t$, we have
  $(Q_r'+s)\cap(Q_r'+t) = \emptyset$. To see this, suppose that $s,t\in Q_{q/r}'$ and that $(Q_r'+s)\cap(Q_r'+t)$ is nonempty. These two conditions imply that $\abs{s-t} = b/(q/r)$ for some $b=0, \dots, q/r-1$, and that $\abs{s-t} = a/r$ for some $a = 0, \dots, r-1$. Since $\gcd(r,q/r)=1$, it must be that $a=b=0$, hence, $t=s$.
  
  This shows that
  \begin{equation*}
    \#(Q_r' + Q_{q/r}') = (\#Q_r')\cdot(\#Q_{q/r}') = \#Q_q'
  \end{equation*}
 and the lemma follows.
\end{proof}

\begin{lem}\label{lem:eps}
  For every $\eps>0$ and any $Q>0$, there exist distinct primes
  $Q < p_1 < \dots < p_\ell$ such that
  \begin{equation*}
    \frac{\varphi(P)}{P} < \eps,
  \end{equation*}
  where $P = \prod_{i=1}^\ell p_i$.
\end{lem}

\begin{proof}
  Let $Q < p_1 < \dots < p_\ell$ be consecutive primes. Note that
  \begin{equation*}
    \frac{\varphi(P)}{P} = \prod_{i=1}^\ell \frac{p_i-1}{p_i} =     \prod_{i=1}^\ell \parens*{1 - \frac{1}{p_i}}.
  \end{equation*}
  Since $\sum_{p \in \mathbb{P}} \frac{1}{p} = \infty$, choosing $\ell$ sufficiently large proves the claim.
\end{proof}

\begin{proof}[Proof of Theorem \ref{thm:weakcounterex}]

The function $\psi:\NN\to[0,\infty)$ will be supported in blocks of
integers:
\begin{equation*}
  \supp \psi = B_1 \cup B_2 \cup \dots\cup B_j\cup\cdots \subset \NN.
\end{equation*}
To start the construction, put $P_0 = 1$.

For each $j\geq 1$, by Lemma \ref{lem:eps} we may choose primes
$P_{j-1}< p_1 < p_2 < \dots < p_\ell$ such that
\begin{equation}\label{eq:smallmeas}
 \frac{\varphi(P_j)}{P_j} \leq 2^{-j},
\end{equation}
where $P_j = \prod_{i=1}^\ell p_i$.  Defining
\begin{equation*}
  B_j = \set{q \in\NN : q\mid P_j, \; q>1},
\end{equation*}
we observe that the blocks $B_j$ are disjoint from one another.

We now define $(y_q)_{q\in\NN}$ and $\psi:\NN\to [0,\infty)$. For
each $q\in B_j$, let $y_q\in\RR$ be such that
\begin{equation*}
 \frac{y_q}{q} \in Q_{P_j/q}'.
\end{equation*}
(The values of $y_q$ for $q\notin \bigcup_j B_j$ may remain
arbitrary.) Define
\begin{equation*}
  \psi(q)
  =
  \begin{cases}
    \frac{q}{2P_j} &\textrm{if } q \in B_j,\\
    0 &\textrm{otherwise.}
  \end{cases}
\end{equation*}
We claim that
\begin{equation}
  \label{eq:sumdiv}
  \sum_{q=1}^\infty \frac{\varphi(q)\psi(q)}{q} = \infty
\end{equation}
and that $\lambda(W'(y,\psi))=0$.\\

The first claim follows from calculating
\begin{align*}
  \sum_{q=1}^\infty \frac{\varphi(q)\psi(q)}{q}
  &= \sum_{j=1}^\infty \sum_{q\in B_j}\frac{\varphi(q)\psi(q)}{q} 
  = \sum_{j=1}^\infty \frac{1}{2P_j}\sum_{\substack{q \mid P_j \\ q > 1}}\varphi(q) 
  = \sum_{j=1}^\infty \frac{P_j - 1}{2P_j} = \infty,
\end{align*}
so we are left to show that $\lambda(W'(y,\psi))=0$. Notice that
\begin{equation}\label{eq:limsupj}
  W'(y,\psi) = \limsup_{q\to\infty} \Aq = \limsup_{j\to\infty}\parens*{\bigcup_{q\in B_j}A_q^{y}}.
\end{equation}
Furthermore, observe that for each $q\in B_j$, we have by Lemma \ref{lem:sumset}
\begin{equation*}
 \Aq \subset Q_{P_j}' + \parens*{-\frac{1}{2P_j},\frac{1}{2P_j}}.
\end{equation*}
Therefore by~\eqref{eq:smallmeas},
\begin{equation}\label{eq:3}
  \lambda\parens*{\bigcup_{q \in B_j} \Aq} \leq \frac{\varphi(P_j)}{P_j} \leq 2^{-j}.
\end{equation}
 Since $\sum_{j\in \N} 2^{-j}$ converges, the convergence Borel--Cantelli Lemma gives 
\begin{equation*}
 \lambda\parens*{\limsup_{j\to\infty}\parens*{\bigcup_{q\in B_j}\Aq}}=0,
\end{equation*}
which by~(\ref{eq:limsupj}) implies $\lambda(W'(y,\psi))=0$, finishing
the proof.
\end{proof}

\begin{rem}
    We note that in the specific instance of $\psi$ chosen above, we do not contradict the possibility of Question \ref{conj:weakdscmoving} being true: Since $\psi(P_j) = \frac{1}{2}$ for every $j \in \N$, this trivially implies that $W(y,\psi)$ has full measure. Note that the value $\frac{1}{2}$ is chosen arbitrarily and can be replaced by any positive constant $c > 0$, in which case we would still obtain by a well-known fact (for a modern formulation see e.g. \cite[Theorem 4]{BHV24}), $\lambda(W(y,\psi)) =1$.
    \end{rem}

\section{Proof of Theorem \ref{thm:bootstrapped}}
\label{sec:bootstrap}

Let $\y = (\byq)_{q \in \mathbb{N}}$ be as in the statement of Theorem \ref{thm:bootstrapped}. The convergence case of the theorem is standard, so let $\psi:\ZZ^n\to[0,\infty)$ satisfy
\begin{equation*}   \sum_{\q\in\ZZ^n\setminus\set{0}}\parens*{\frac{\varphi(\gcd(\q))\psi(\q)}{\gcd(\q)}}^m = \infty.
\end{equation*}
Clearly, there is some orthant of $\ZZ^n$ to which we may restrict the sum without losing the divergence, and without loss of generality we may suppose that it is the positive orthant, so that
\begin{equation}\label{eq:bootsum}
\sum_{\q\in\ZZ_{\geq 0}^n\setminus\set{0}}\parens*{\frac{\varphi(\gcd(\q))\psi(\q)}{\gcd(\q)}}^m = \infty.
\end{equation}
For $Q\geq 1$, define $\Psi:\NN\to [0,\infty]$ by
\begin{equation}
    \Psi_Q(d) = \Bigg(\sum_{\substack{\q\in\ZZ_{\geq 0}^n \\ \gcd(\q)=d \\ \abs{\q/d}\geq Q}}\psi(\q)^m\Bigg)^{1/m}.
\end{equation}
By~\cite[Proposition~2.1]{ramirez2024metric}, we have $\lambda(W_{n,m}'(\y,\psi))=1$ as soon as there exists some $d,Q$ with $\Psi_Q(d)=\infty$, so we may assume that $\Psi_Q(d)$ is always finite.

Next, suppose that there exists some $\q\in\ZZ_{\geq 0}^n$ with $\gcd(\q)=1$ such that
\begin{equation}\label{eq:raysum} 
\sum_{d = 1}^\infty\parens*{\frac{\varphi(d)\psi(d\q)}{d}}^m = \infty.
\end{equation}
Then the function $\psi_\q(d):=\psi(d\q)$ satisfies the divergence condition in Theorem \ref{thm:movingdsc}, and therefore, $\lambda(W_m'(\y,\psi_\q))=1$. But
\begin{equation*}
    T^{-1}\parens{W_m'(\y,\psi_\q)} \subset W_{n,m}'(\y,\psi),
\end{equation*}
where $T:[0,1)^{nm}\to[0,1)^m$ by $\X\mapsto\q\X \pmod 1$. This mapping is measure preserving, therefore $\lambda(W_{n,m}'(\y,\psi))=1$ and we are done. 

Thus we are left to consider the case where there is no $\q\in\ZZ_{\geq 0}^n$ for which~\eqref{eq:raysum} holds. Let $Q\geq 1$. Then  
\begin{align*}
    \sum_{d=1}^\infty \parens*{\frac{\varphi(d)\Psi_Q(d)}{d}}^m 
    &= \sum_{d=1}^\infty \parens*{\frac{\varphi(d)\Psi_1(d)}{d}}^m - \sum_{d=1}^\infty \sum_{\substack{\q\in\ZZ_{\geq 0}^n \\ \gcd(\q)=d \\ \abs{\q/d} < Q}}\parens*{\frac{\varphi(d)\psi(\q)}{d}}^m \\
    &= \sum_{d=1}^\infty \parens*{\frac{\varphi(d)\Psi_1(d)}{d}}^m - \sum_{\substack{\q\in\ZZ_{\geq 0}^n \\ \gcd(\q)=1 \\ \abs{\q} < Q}} \sum_{d=1}^\infty \parens*{\frac{\varphi(d)\psi(d\q)}{d}}^m
\end{align*}
and the assumption implies that the second term is finite. Meanwhile,~\eqref{eq:bootsum} is equivalent to
\begin{align*}
    \sum_{d=1}^\infty \parens*{\frac{\varphi(d)\Psi_1(d)}{d}}^m = \infty,
\end{align*}
so $\Psi_Q$ satisfies the divergence condition in Theorem \ref{thm:movingdsc}, which implies that $\lambda(W_m'(\y,\Psi_Q))=1$. 
Since $Q>1$ was arbitrary, we can use~\cite[Theorem~1.2]{ramirez2024metric}, which states that
\begin{equation}\label{eq:bootstrap}
    \inf_{Q} \lambda\left(W_m'(\y,\Psi_Q)\right)>0 \qquad\implies\qquad \lambda\left(W_{n,m}'(\y,\psi)\right)=1, 
\end{equation}
and this proves Theorem \ref{thm:bootstrapped}.\qed

\appendix

\section{Proof of Theorem \ref{thm:movingtarget}}
\label{sec:two}

As remarked in Section~\ref{sec:remarks}, the $m\geq 3$ cases of Theorem \ref{thm:movingtarget} are found in~\cite{allen2021independence,Yu}, but we know of no reference that mentions the two-dimensional Khintchine Theorem in the context of moving targets. Here we provide a short proof that works for all $m\geq 2$. It would be possible to give a proof that uses the overlap estimates Lemma \ref{lemma_over}, but we present a proof that uses instead classical estimates that are sometimes referred to as the trivial overlap estimates.

\begin{lem}[Trivial overlap estimates]
  For $1\leq r < q$ and $y,z\in\RR$,
  \begin{equation*}
        \lambda_1(A_q^y\cap A_r^z) \ll \psi(q)\psi(r) + \frac{\psi(q)}{q}\varphi(\gcd(q,r)). 
  \end{equation*}
\end{lem}

\begin{proof}[Sketch of proof]
  The idea is that an intersection between an interval from $A_q^y$
  and an interval from $A_r^z$ can have measure at most $\delta$, where
  \begin{equation*}
    \delta = 2\min\parens*{\frac{\psi(q)}{q}, \frac{\psi(r)}{r}}.
  \end{equation*}
  In order to bound the measure of $A_q^y\cap A_r^z$, we estimate the
  number of intersections between intervals of $A_q^y$ and intervals
  of $A_r^z$, and multiply by $\delta$. Such intersections occur
  precisely when the centers of the intervals have a distance between
  them less than $\tfrac{\Delta}{2} + \tfrac{\delta}{2}$, where
  \begin{equation*}
    \Delta = 2\max\parens*{\frac{\psi(q)}{q}, \frac{\psi(r)}{r}}.
  \end{equation*}
  Relaxing a bit, we seek a bound on
  \begin{multline*}
    \set*{(a,b) : 1\leq a \leq q, 1\leq b \leq r, (a,q)=1, (b,r)=1, \abs*{\frac{a+y}{q} - \frac{b+z}{r}} \leq \Delta} \\
    = \set*{(a,b) : 1\leq a \leq q, 1\leq b \leq r, (a,q)=1, (b,r)=1, ar-bq \in I}     
  \end{multline*}
  where $I$ is an interval of length $2\Delta q r$.  There are
  $\leq 2\Delta qr / \gcd(q,r) + 1$ integers $c=ar-bq$ that can appear
  in $I$, and each one has $\leq \varphi(\gcd(q,r))$ representations
  in the form $ar-bq=c$ with $(a,q)=1$ and $(b,r)=1$. Putting things
  together, one obtains
  \begin{align*}
    \lambda_1(A_q^y\cap A_r^z)
    &\ll  \delta \parens*{\frac{2\Delta q r}{\gcd(q,r)} + 1} \varphi(\gcd(q,r)) \\
    &\ll \psi(q)\psi(r) + \frac{\psi(q)}{q}\varphi(\gcd(q,r)),
  \end{align*}
  as needed.
\end{proof}

A corollary analogous to Corollary \ref{multidim_statements} follows.
\begin{cor}
    For $1\leq r < q$ and a sequence $\y = (\y_q)_{q\in\NN}\subset \RR^m$, 
  \begin{equation}
    \label{eq:trivial}
    \lambda_m(\Aqm\cap \Arm) \ll_m \psi(q)^m\psi(r)^m + \parens*{\frac{\psi(q)}{q}}^m\varphi(\gcd(q,r))^m. 
  \end{equation}
\end{cor}

\begin{proof}
The proof is identical to the proof of~\eqref{multidim_overlap}.
\end{proof}

\begin{proof}[Proof of Theorem \ref{thm:movingtarget}]
Let $\y = (\y_q)_{q \in \N}\subset \RR^m$ and $\psi:\NN\to [0,\infty)$
a monotonic function such that $\sum \psi(q)^m$ diverges. Assume that $\psi(q)^m \to 0$, for otherwise the
result is trivial. With that assumption, no generality is lost in
further assuming that $\psi(q) \leq 1/2$ for all $q$, so that we have
\begin{equation}\label{eq:trivmeas}
    \lambda_m(\Aqm) = \parens*{\frac{2\varphi(q)\psi(q)}{q}}^m
\end{equation}
for all $q\in\NN$.

Notice that $W_m(\y,\psi)$ contains $\limsup_{q\to\infty}\Aqm$ (which we recall is the set of \textit{coprime} approximations). It is therefore enough to show that 
$
    \lambda_m(\limsup_{q \to \infty} \Aqm) =1.
$
From~\eqref{eq:trivial} we have
\begin{equation}
    \lambda_m(\Aqm\cap \Arm) \ll \psi(q)^m\psi(r)^m + \frac{\psi(q)^m}{q^m}\varphi(\gcd(q,r))^m.
\end{equation}
Now, for $Q>1$ we have 
  \begin{align*}
    \sum_{q,r = 1}^Q\lambda(\Aqm\cap \Arm)
    &\ll \sum_{q,r=1}^Q \psi(q)^m\psi(r)^m + \sum_{q=1}^Q\frac{\psi(q)^m}{q^m} \sum_{r=1}^q\varphi(\gcd(q,r))^m \\
    &\ll \parens*{\sum_{q=1}^Q \psi(q)^m}^2 + \sum_{q=1}^Q \psi(q)^m,
  \end{align*}
  by Lemma \ref{lem:phigcd}. Since $\psi(q)$ is nonincreasing and
  $\sum_{n \leq N} \frac{\varphi(n)}{n} \gg N$, summation by parts yields
  \begin{equation*}
      \sum_{q=1}^Q \psi(q)^m \ll \sum_{q=1}^Q \parens*{\frac{\varphi(q)\psi(q)}{q}}^m,
  \end{equation*}
  so we have established
  \begin{align*}
    \sum_{q,r = 1}^Q\lambda_m(\Aqm \cap \Arm) \ll \parens*{\sum_{q=1}^Q \parens*{\frac{\varphi(q)\psi(q)}{q}}^m}^2 \overset{\eqref{eq:trivmeas}}{\ll} \parens*{\sum_{q=1}^Q \lambda_m(\Aqm)}^2.
  \end{align*} 
  Again, the theorem follows from Propositions \ref{BC5} and \ref{prop:mixing}.
\end{proof}

\section{Proof of Lemma \ref{pv_estimate} for $m \geq 2$}\label{sec:PVredo}

This section is a recreation of arguments from Pollington--Vaughan~\cite{PV90}, included in this paper for completeness and convenience.

\begin{lem}\label{lem:erdos}
For given $s \in \mathbb{N}$, let 
    \[g(s) := \min\Bigg\{n \in \N: \sum_{\substack{p \mid s\\p > n}} \frac{1}{p} < \frac{1}{2}\Bigg\}.\]
Then we have uniformly in $x \geq 1, v \geq 1$,
\[\#\left\{n < x: g(n) = v\right\}
 \ll \frac{x}{(v-1)!}.
\]
\end{lem}

\begin{proof}
Note that by the definition of $g$, if $g(b) = v$, we obtain that
$\sum\limits_{\substack{p \mid n\\ p > v-1}} \frac{1}{p} \geq \frac{1}{2}$.
By~\cite[Lemma 3]{E70}, we have 
\[\#\Bigg\{n \leq x: \sum_{\substack{p \mid n\\ p > v-1}} \frac{1}{p} \geq \frac{1}{2}\Bigg\}
\ll \frac{x}{(v-1)!},
\]
which proves the statement.
\end{proof}

\begin{proof}[Proof of Lemma \ref{pv_estimate} for $m \geq 2$]
In the following, we will always understand $\log(x) := \max\{1, \Log(x)\}$ where $\Log$ denotes the natural logarithm. With $g$ as in Lemma \ref{lem:erdos}, Mertens' Theorem shows that
\begin{equation}\label{bound_on_1p}
    \sum_{\substack{p \mid s}} \frac{1}{p} \leq \frac{1}{2} + \sum_{p \leq g(s)} \frac{1}{p} 
    \ll \log \log g(s).
\end{equation}
    Thus writing
    \[t= t(q,r) := \max\left\{g\big(\tfrac{q}{(q,r)}\big),g\big(\tfrac{r}{(q,r)}\big)\right\},\]
    we observe that

\[
\prod_{\substack{p \mid \frac{qr}{(q,r)^2}\\ p > D(q,r)}}\left(1 + \frac{1}{p}\right)
\leq \exp\Bigg(
\sum_{\substack{p \mid \frac{q}{(q,r)}\\p \leq D(q,r)}} \frac{1}{p} + \sum_{\substack{p \mid \frac{r}{(q,r)}\\p > D(q,r)}} \frac{1}{p} \Bigg) \ll \begin{cases} 1 &\text{ if } t \leq D(q,r),\\
    \log t, &\text{ if } t \geq D(q,r).
\end{cases}
\]
This implies that

\begin{equation}\label{casesPV}
    M(q,r)^m \ll \begin{cases}
        \productm &\text{ if } t \leq D(q,r),\\
        \productm &\text{ if } D(q,r) \leq 1,\\
       \productm(\log t)^m &\text{ if } 1 < D(q,r) < t.
    \end{cases}
\end{equation}
Note that summing over all pairs $(q,r)$ that belong to the first two cases of \eqref{casesPV} yields a contribution $\ll \left(\sum_{q \leq Q} \left(\frac{\varphi(q)\psi(q)}{q}\right)^m\right)^2$, and therefore, we are only left to consider pairs from the case where $1 < D(q,r) < t$. By symmetry reasons, it suffices to consider pairs where $r < q$.\\

For the remainder of the proof, we let $a = \frac{q}{(q,r)}$, $b = \frac{r}{(q,r)}$ and recall
$\Delta =  \max\left\{\frac{\psi(q)}{q},\frac{\psi(r)}{r}\right\}, D(q,r) = \frac{\Delta qr}{(q,r)}$.
Note that if $\frac{\psi(q)}{q} > \frac{\psi(r)}{r}$, then 
\begin{equation}
    \label{bigq}
D(q,r) = b\psi(q) > a\psi(r)
\end{equation}
and if 
$\frac{\psi(q)}{q} \leq \frac{\psi(r)}{r}$, then 
\begin{equation}\label{bigr}D(q,r) = a\psi(r) \geq b\psi(q).\end{equation}

Note that for $r < q$, it is impossible to satisfy 
$\frac{\psi(r)}{r} < \frac{\psi(q)}{q}$ and $\psi(r) > \psi(q)$ simultaneously.
Thus the following case distinction is exhaustive:
\[
\begin{split}
\mathcal{E}_1 &:= \left\{(q,r) \in [1,Q]^2:\quad  r < q,\quad \Delta = \frac{\psi(r)}{r},\quad 1 < a\psi(r) < t,\quad \psi(r) \leq \psi(q) \right\},\\
\mathcal{E}_2 &:= \left\{(q,r) \in [1,Q]^2:\quad  r < q,\quad \Delta = \frac{\psi(r)}{r},\quad
 1 < a\psi(r) < t,\quad
\psi(r) > \psi(q) \right\},\\
\mathcal{E}_3 &:= \left\{(q,r) \in [1,Q]^2:\quad  r < q,\quad \Delta = \frac{\psi(q)}{q},\quad 1 < b\psi(q) < t,\quad \psi(r) \leq \psi(q) \right\}.\\
\end{split}
\]



\subsection*{Contributions from $\mathcal E_1$}  

We decompose
\begin{multline*}
\sum_{(q,r) \in \mathcal{E}_1} M(q,r)^m \ll \sum_{(q,r) \in \mathcal{E}_1} \measqm \measrm (\log t)^m \\
\ll
\sum_{t \in \N} (\log t)^{m}
\sum_{\substack{(u,v) \in \N\\ \max\{u,v\} = t}}
\sum_{\substack{(q,r) \in \mathcal{E}_1,\\ g(a) = u,\\ g(b) = v}} \measqm \measrm.
\end{multline*}
By $a\psi(r) < t$, we have
\begin{equation}\label{case11}\measr \leq \psi(r) \leq \frac{t}{a}.\end{equation}
Furthermore, using $\psi(r) \leq \psi(q)$ and $a \psi(r) > 1$, we have 
\begin{equation}\label{case12}a > \frac{1}{\psi(r)} \geq \frac{1}{\psi(q)}.\end{equation}
Finally, using \eqref{bigr} and $a\psi(r) < t$, we get
\begin{equation}\label{case13}b \leq \frac{a\psi(r)}{\psi(q)} < \frac{t}{\psi(q)}.
\end{equation}
Note that for $q,a,b$ fixed, there exists at most one $r$ such that 
\[a = \frac{q}{\gcd(q,r)},\quad b = \frac{r}{\gcd(q,r)}.\]
For any fixed pair $(u,v)$ with $t = \max\{u,v\}$, we therefore get

\begin{equation}\label{suvbound}\begin{split}&(\log t)^m\sum_{\substack{(q,r) \in \mathcal{E}_1\\ g(a) = u,\\ g(b) = v}} \measqm \measrm \\\leq &\;(\log t)^m \sum_{q \leq Q} \measqm \sum_{\substack{b < \frac{t}{\psi(q)}\\ g(b) = v}}
\sum_{\substack{a < \frac{t}{\psi(q)}\\ g(a) = u}} \measrm
\\\leq &\;(\log t)^m \sum_{q \leq Q} \measqm 
\sum_{\substack{a > \frac{1}{\psi(q)}\\ g(a) = u}}
\sum_{\substack{b < \frac{t}{\psi(q)}\\ g(b) = v}}
 \left(\tfrac{t}{a}\right)^m
\\=&\; (t \log t)^m \sum_{q \leq Q} \measqm 
\sum_{\substack{a > \frac{1}{\psi(q)}\\ g(a) = u}} \frac{1}{a^m}
\sum_{\substack{b < \frac{t}{\psi(q)}\\ g(b) = v}}1.
\end{split}
\end{equation}
Using again Lemma \ref{lem:erdos}, employing summation by parts, we obtain
\begin{equation}\begin{split}
\sum_{\substack{a > \frac{1}{\psi(q)}\\g(a) = u}} \frac{1}{a^m}
&\leq \sum_{t > \frac{1}{\psi(q)}}\left(\sum_{1< j < t} \mathds{1}_{[g(j) = u]}\right)
\left(\frac{1}{t^m} - \frac{1}{(t+1)^{m}}\right)
\\&\ll \sum_{t > \frac{1}{\psi(q)}} \frac{t}{(u-1)!} \frac{1}{t^{m+1}}
\\&\ll \frac{\psi(q)^{m-1}}{(u-1)!}.
\end{split}
\label{anatomy2}
\end{equation}
Combining \eqref{suvbound}, \eqref{anatomy2} and another application of Lemma \ref{lem:erdos} for the $b$-summation yields

\begin{equation}
\begin{split}
    &\;(\log t)^m\sum_{\substack{(q,r) \in \mathcal{E}_1\\ g(a) = u,\\ g(b) = v}} \measqm \measrm 
\sum_{\substack{b < \frac{t}{\psi(q)}\\ g(b) = v}}1
\\\ll&\; (t \log t)^m \sum_{q \leq Q} \measqm \frac{\frac{t}{\psi(q)}}{(v-1)!}\frac{\psi(q)^{m-1}}{(u-1)!}
\\\leq&\; \frac{(u+v)^2 (\log (u+v))^m}{(u-1)!(v-1)!} \sum_{q \leq Q} \measqm \psi(q)^{m-2}.
\end{split}
\end{equation}
Since \[\sum_{u,v \geq 1} \frac{(u+v)^2 (\log (u+v))^m}{(u-1)!(v-1)!} \ll 1,
\]
$m \geq 2$ and $\psi(q) \leq 1$, we obtain 
\begin{equation}\label{eq:bing}
\sum_{(q,r) \in \mathcal{E}_1}M(q,r)^m
\ll \sum_{q \leq Q} \measqm.
\end{equation}

\subsection*{Contributions from $\mathcal E_2$}
Recall
\[\mathcal{E}_2 := \left\{(q,r) \in [1,Q]^2:\quad  r < q,\quad \Delta = \frac{\psi(r)}{r},\quad
 1 < a\psi(r) < t,\quad
\psi(r) > \psi(q) \right\}.\]
By $a\psi(r) < t$ and $\psi(r) > \psi(q)$, we have
\begin{equation}\label{case21}a < \frac{t}{\psi(r)} < \frac{t}{\psi(q)}.
\end{equation}
as well as by $r < q$,
\begin{equation}\label{case22}
    b = \frac{r}{(q,r)} < \frac{q}{(q,r)} < \frac{t}{\psi(q)}.
\end{equation}
Therefore, reparametrizing the pair $(q,r)$ by $(r,a,b)$ and using $\psi(q) < \psi(r)$, we get
\begin{equation}\begin{split}
    &\;(\log t)^m\sum_{\substack{(q,r) \in \mathcal{E}_2\\ g(a) = u,\\ g(b) = v}} \measqm \measrm
\\\leq &\;(\log t)^m\sum_{r \leq Q} \measrm \sum_{\substack{a < \frac{t}{\psi(q)}\\g(a) = u}}\sum_{\substack{b < \frac{t}{\psi(q)}\\g(b) = v}}\psi(r)^m
\\= &\;(\log t)^m\sum_{r \leq Q} \measrm \psi(r)^m 
\#\left\{n < \frac{t}{\psi(q)}: g(n) = u\right\} \cdot \#\left\{n < \frac{t}{\psi(q)}: g(n) = v\right\}
\\\ll&\; \frac{t^2(\log t)^m}{(u-1)!(v-1)!}\sum_{r \leq Q} \measrm \psi(r)^{m-2},
\end{split}
\end{equation}
where we used Lemma \ref{lem:erdos} twice. Summing over $u,v$ again yields 
\begin{equation}\label{eq:bang}
\sum_{(q,r) \in \mathcal{E}_2}M(q,r)^m
\ll \sum_{q \leq Q} \measqm.
\end{equation}

\subsection*{Contributions from $\mathcal E_3$}
Recall
\[\mathcal{E}_3 := \left\{(q,r) \in [1,Q]^2:\quad  r < q,\quad \Delta = \frac{\psi(q)}{q},\quad 1 < b\psi(q) < t,\quad \psi(r) \leq \psi(q) \right\}.\]
 By $1 < b\psi(q) < t$ and $r< q$, we get
 \[a = \frac{q}{(q,r)} > \frac{r}{(q,r)} = b > \frac{1}{\psi(q)}, \quad b < \frac{t}{\psi(q)}.\]
Furthermore, by $\frac{\psi(r)}{r} \leq \frac{\psi(q)}{q}$, we have
\begin{equation*}
    \measr \leq \psi(r) \leq \frac{\psi(q)r}{q} = \frac{b}{a}\psi(q).
\end{equation*}
This leads to
\[
\begin{split}
&\;(\log t)^m\sum_{\substack{(q,r) \in \mathcal{E}_1\\ g(a) = u,\\ g(b) = v}} \measqm \measrm \\\ll& \;(\log t)^m\sum_{q \leq Q} \measqm  \sum_{\substack{a < \frac{1}{\psi(q)}\\g(a) = u}}\sum_{\substack{b < \frac{t}{\psi(q)}\\g(b) = v}} \left(\frac{b}{a}\psi(q)\right)^m
\\=&\; (\log t)^m\sum_{q \leq Q} \measqm \psi(q)^m \Bigg(\sum_{\substack{a < \frac{1}{\psi(q)}\\g(a) = u}}\frac{1}{a^m}\Bigg)\Bigg(\sum_{\substack{b < \frac{t}{\psi(q)}\\g(b) = v}}b^m\Bigg).
\end{split}
\]
Another summation by parts combined with Lemma \ref{lem:erdos} gives 
\[\sum_{\substack{b < \frac{t}{\psi(q)}\\g(b) = v}}b^m
\leq \sum_{\substack{j < \frac{t}{\psi(q)}}} \left(\sum_{b \leq j+1} \mathds{1}_{[g(b) = v]}\right)\left((j+1)^m - j^m\right) \ll \frac{1}{(v-1)!}\sum_{j < \frac{t}{\psi(q)}} j^m \ll \frac{t^{m+1}}{\psi(q)^{m+1}(v-1)!},
\]
and using \eqref{anatomy2}, we deduce
\[(\log t)^m\sum_{\substack{(q,r) \in \mathcal{E}_1\\ g(a) = u,\\ g(b) = v}} \measqm \measrm \ll \frac{t^{m+1}(\log t)^m}{(u-1)!(v-1)!}\sum_{q \leq Q} \measqm \psi(q)^{m-2}.\]
Summation over $u,v$ again yields
\begin{equation}\label{eq:boom}
\sum_{(q,r) \in \mathcal{E}_3} M(q,r)^m \ll \sum_{q \leq Q} \measqm
\end{equation}
and the proof of the lemma follows from combining~\eqref{eq:bing},~\eqref{eq:bang}, and~\eqref{eq:boom} with the contributions from the first two cases of~\eqref{casesPV}.
\end{proof}

\subsection*{Acknowledgements} 
MH was supported by the EPSRC grant EP/X030784/1. FAR thanks the Number Theory Group at the University of York for their hospitality.

\bibliographystyle{plain}
\bibliography{bibliography.bib}

\end{document}